# Applications of the Digital-Discrete Method in Smooth-Continuous Data Reconstruction

Li Chen
Department of Computer Science and Information Technology
University of the District of Columbia
Email: *lchen@udc.edu*

**Abstract:** This paper presents some applications using recently developed algorithms for smooth-continuous data reconstruction based on the digital-discrete method. The classical discrete method for data reconstruction is based on domain decomposition according to guiding (or sample) points. Then the Spline method (for polynomial) or finite elements method (for PDE) is used to fit the data. Our method is based on the gradually varied function that does not assume the property of being linearly separable among guiding points, i.e. no domain decomposition methods are needed. We also demonstrate the flexibility of the new method and its potential to solve a variety of problems. The examples include some real data from water well logs and harmonic functions on closed 2D manifolds. This paper presents the results from six different algorithms. This method can be easily extended to higher multi-dimensions. We also include an advanced consideration related to the use of gradually varied mapping.

## 1. Introduction

Data reconstruction is used to fit a function based on the observations of some sample (guiding) points. This paper presents some applications of using recently developed algorithms for smooth-continuous data reconstruction based on the digital-discrete method [2]. The classical discrete method for data reconstruction is based on domain decomposition according to guiding (or sample) points. Then, the Spline method (for polynomial) or finite elements method (for PDE) is used to fit the data.

Some successful methods have been discovered or proposed to solve the problem including the Voronoi-based surface method and the moving least square method [1][11][14][16][19].

Our method is based on the gradually varied function that does not assume the property of the linearly separable among guiding points, i.e. no domain decomposition methods are needed. We also demonstrate the flexibility of the new method and its potential to solve a variety of problems. The examples include some real data from water well logs and harmonic functions on closed 2D manifolds. This paper presents the results from six different algorithms. This method can be easily extended to higher multi-dimensions.

This paper will mainly present applications of using the method in many different cases including rectangle domains and closed surfaces. We directly deal with smooth functions on the (same) piecewise linear manifold or a non-Jordon graph. We have also applied this method to groundwater flow equations [4]. We discovered that the major difference between our methods and existing methods is that the former is a true nonlinear approach [3]. We analyzed the relationships and differences among the different methods in [3].

## 2. Basic Concepts

This this section, we first introduce the basic concepts and then provide some information about the existing theories. Gradual variation is a discrete method that can be built on any graph [7-8]. The gradually varied surface is a special discrete surface.

Let $A_1 < A_2 < \ldots < A_n$. The Concept of Gradual Variation: Let function f: D$\rightarrow${$A_1, A_2, \ldots, A_n$}. If a and b are adjacent in D, then it is implied that f(a)=f(b), or f(b) =A(i-1) or A(i+1) when f(a)=Ai. Point (a,f(a)) and (b,f(b)) are then said to be gradually varied. A 2D function (surface) is said to be gradually varied if every adjacent pair is gradually varied.

Discrete Surface Fitting: Given J$\subseteq$D, and f: J$\rightarrow${$A_1,A_2,\ldots A_n$}, decide if there is exists an F: D$\rightarrow${$A_1,A_2,\ldots,A_n$} such that F is gradually varied where f(x)=F(x), x in J.

**Theorem** (Chen, 1989) [7-8] The necessary and sufficient conditions for the existence of a gradually varied extension F is: for all x,y in J, d(x,y)$\geq$ |i-j|, f(x)=Ai and f(y)=Aj, where d is the distance between x and y in D.

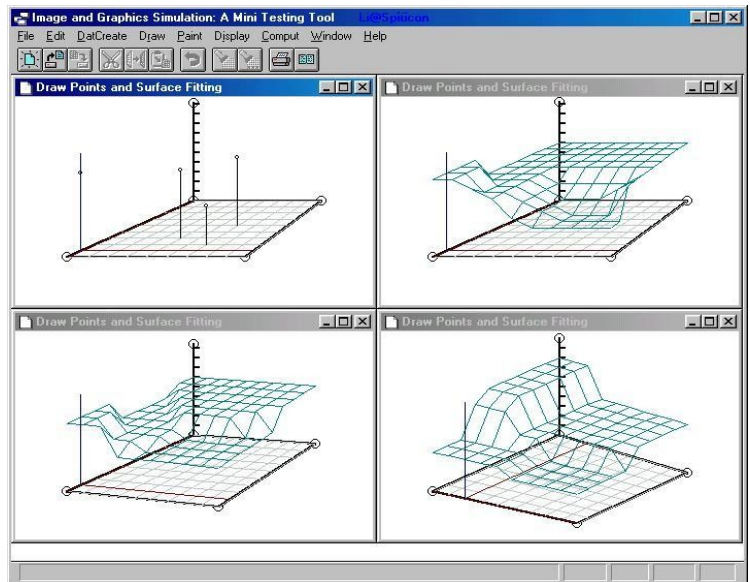
Fig. 2 Examples of gradually varied functions

A gradually varied surface fitting software component was included in a lab-use oriented software system in 1997. The software could demonstrate the arbitrary guiding points of a gradually varied surface fitting in a 10x10 grid domain.

In 2004, Chen and Adjei proposed a method to do continuous and differential surface fitting based on lambda-connectedness that was a continuous space treatment of discrete functions [9] [5]. However, this method has not been implemented. In 2008, Chen implemented an algorithm for groundwater data reconstruction; the implementation was not fully successful. The reason was that the Lipschitz condition was mainly used for the fitting. In the new algorithm we have used the local Lipschitz function to fit the function [2].

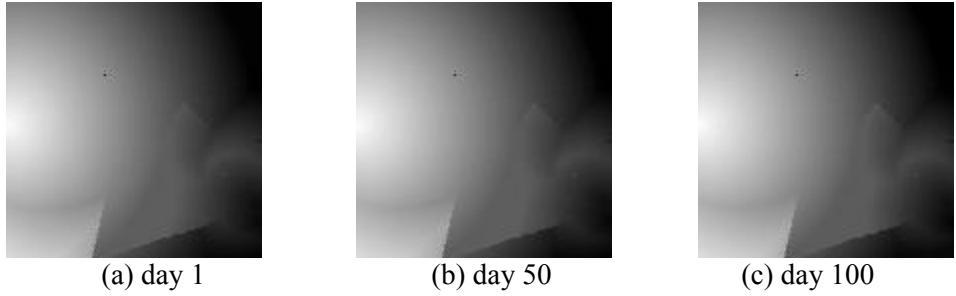

    (a) day 1                (b) day 50             (c) day 100

Fig. 3. Northern VA Groundwater distribution calculated by gradually varied surfaces date from 04/01/07. The intensity indicated the depth of the groundwater. 10 sample points are used.

Thus, the above theorem can be used for a single surface fitting if the condition in the theorem is satisfied. The problem is that the sample data does not satisfy the condition of fitting. So the original algorithm cannot be used directly for individual surface fitting. Another problem is that the theorem is only for "continuous" surfaces. It does not imply a solution for differentiable or smooth functions.

In theory, McShane and Whitney obtained an important theorem for Lipschitz function extension [15][21]. Kirszbraun and later Valentine studied the Lipschitz mapping extension for Hilbert spaces [20]. Recent theoretical research works on these aspects can be found in [12][13][17][18].

The gradually varied function is highly related to the Lipschitz function and the local Lipschitz function. It was proposed for discrete surface and data reconstruction [5]. We will give a brief comparison of the method of gradually varied functions and the McShane-Whitney extension method in Section 6.

## 3. Algorithms and Experiments

In [2], a systematic digital-discrete method for obtaining continuous functions with smoothness to a certain order ($C^n$) from sample data is designed. This method is based on gradually varied functions and the classical finite difference method. The new method has been applied to real groundwater data and the results have validated the method. It is independent from existing popular methods such as the cubic Spline method and the finite element method. The new digital-discrete method has considerable advantages for a large number of real data applications. This digital method also differs from other classical discrete methods that usually use triangulations. The method can potentially be used to obtain smooth functions such as polynomials through its derivatives $f^{(k)}$ along with the solution for partial differential equations such as the harmonic and other important equations.

The new algorithm tries to search for the best solution to the fitting. We have added a component of the classical finite difference method. The major steps of the new algorithm include (This is for 2D functions. For 3D function, we would only need to add a dimension):

Step 1: Load guiding points. In this step we load the data points with observation values.

Step 2: Determine resolution. Locate the points in grid space.
Step 3: Function extension according to the theorem presented in Section 2. This way, we obtain gradually varied or near gradually varied (continuous) functions. In this step the local Lipschitz condition is used.
Step 4: Use finite difference method to calculate partial derivatives. Then obtain the smoothed function.
Step 5: Some multilevel and multi resolution method may be used.

Three sets of real data are tested. The second set uses the same 10 sample points (Fig. 4) and the result is improved. Fig. 5 shows the third data set containing 29 sample points.

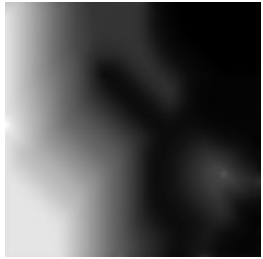

Fig. 4. Northern VA Groundwater distribution calculated by gradually varied surfaces date from 04/01/07. 10 sample points are used at Day 95.

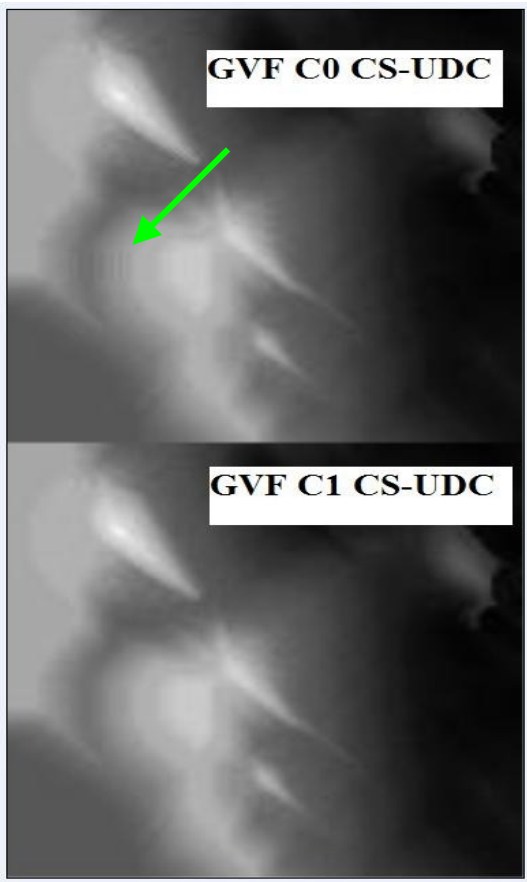

Fig. 5. The picture is the result of fitting based on 29 sample points. The first is a "continuous" surface and the second is the "first derivative." The arrow indicates the interesting area that disappears in the second image.

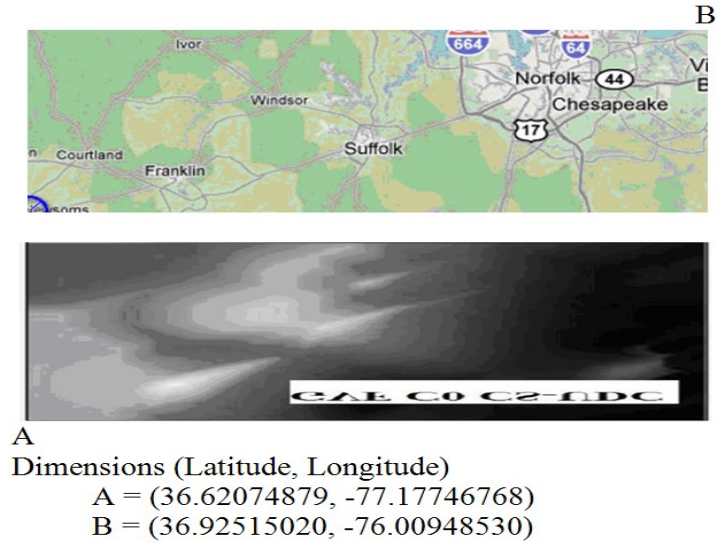

Fig. 6. The map and ground water data

Fig. 6 shows a good match found between the ground water data and the region's geographical map. The brightness of the pixels means the deeper distance from the surface. In mountain areas, the groundwater level is lower in general. Some mismatches may be caused by not having enough sample data points (wells).

## 4. Continuous and Smooth Functions on Manifolds

A gradually varied surface reconstruction does not rely on the shape of domain. And it is not restricted by simplicial decomposition. As long as the domain can be described as a graph, our algorithms will apply. However, the actual implementation will be much more difficult. In the above sections, we have discussed two types of algorithms for a rectangle domain. One is the complete gradually varied function (GVF) fitting and the other is to reconstruct the best fit based on the gradual variation and finite difference method.

The following is the implementation of the method for digital-discrete surface fitting on manifolds (triangulated representation of the domain). The data comes from a modified example in Princeton's 3D Benchmark data sets.

We will have four algorithms related to continuous (and smooth) functions on manifolds. This is because we have 4 cases: (1) ManifoldIntGVF: The GVF extension on point space, corresponding to Delaunay triangulations; the values are integers. (2) ManifoldRealGVF: The GVF fitting on point space, the fitted data are real numbers. (3) ManifoldCellIntGVF: The GVF extension on face (2D-cell) space, corresponding to Voronoi decomposition; the values are integers. (4) ManifoldCellRealGVF: The GVF fitting on face (2D-cell) space, the fitted data are real numbers.

Note that ManifoldRealGVF and ManifoldCellRealGVF are algorithms based on gradually varied functions. The results may or may not be gradually varied. A special data structure was chosen to hold all points, edges, and 2D-cells. It is not very easy to demonstrate the correctness of the reconstruction on manifolds since it is difficult to determine the location on 3D displays.

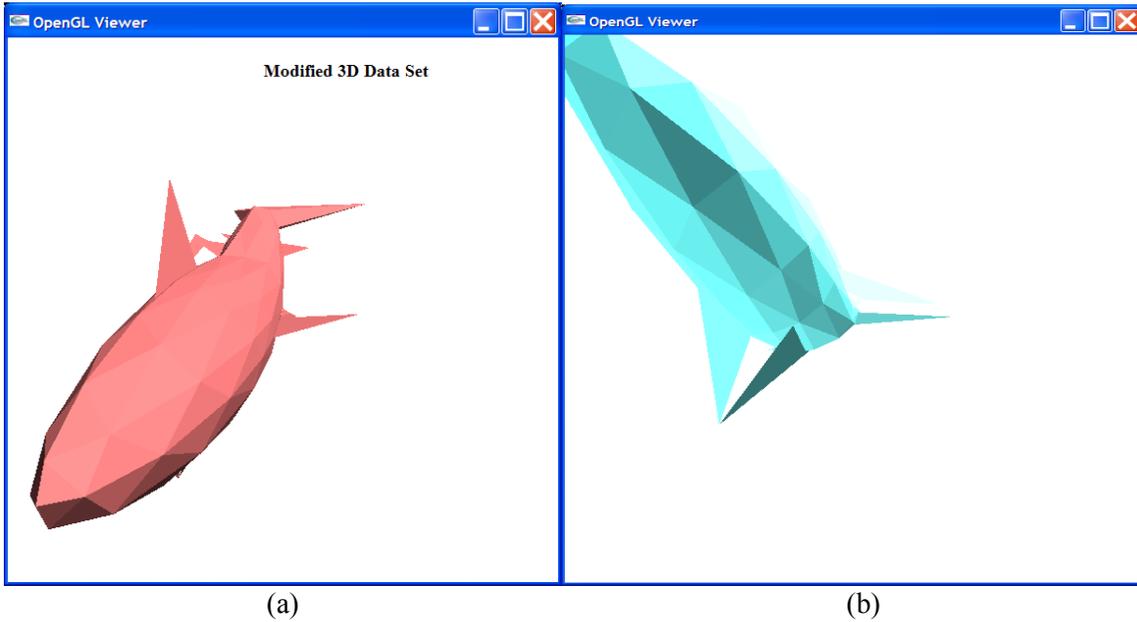

(a) (b)

Fig. 7 The ManifoldIntGVF Algorithm: (a) The left is the original 3D image, (b) the reconstruction based on four sample points.

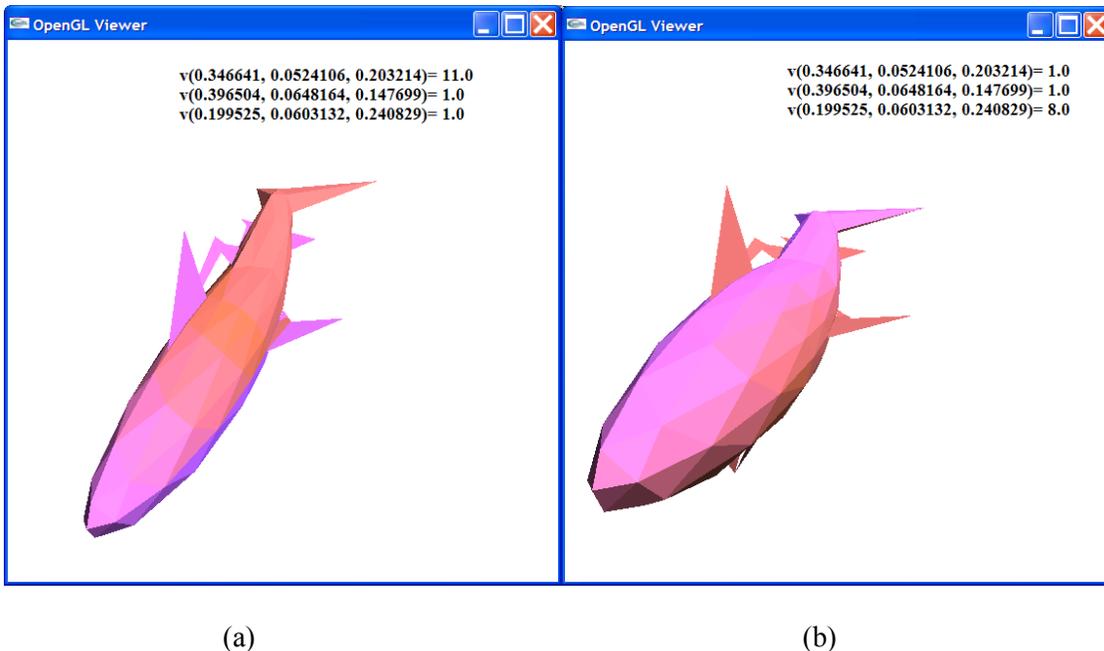

(a) (b)

Fig. 8. The ManifoldIntGVF Algorithm for two 3 sample point sets: (a) The gradual variation is from back to front, (b) The gradual variation from top to the right side.

In Fig. 7 (a) we show an example using ManifoldIntGVF in a 2D closed manifold (can be viewed as the boundary for a 3D object). The image is a fish from the Princeton 3D data sets. We have modified the data by adding some triangles to make the graph connected. Then we put four values on four vertices. Fig. 7 (b) shows the result. In order to avoid being fooled by the 3D display causing the color to change, we have used two additional sets of guiding points. Each of these uses only three guiding points. The second set of guiding points has a pair of values swapped with the first set (The second set of guiding points is just the first set with a pair of values swapped). The values are posted in the pictures. See Fig. 8.

It is very difficult to select values for a relatively large set of sample points for a complete GVF fitting such as ManifoldIntGVF to satisfy the condition of gradual variation (or Lipschitz condition). For ManifoldRealGVF, we can choose any data values. Fig. 9 shows the results with 6 guiding points.

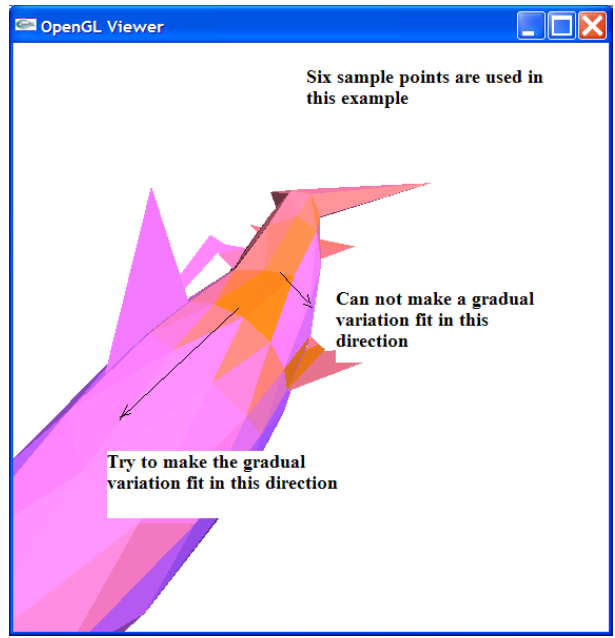

Fig. 9 The ManifoldRealGVF algorithm using 6 guiding points.

Now we can obtain smooth functions based on the GVF fitting. We can use the Catmull-Clark method [11] or other methods. A small problem we encountered is the display. In order to make the display for the vertices, we use the average value of three vertices to be the display value on each triangle. We will make smaller triangles later to make a more precise display. This problem will disappear in the cell-based GVF. However, the graph needs to be calculated first to obtain the cells and their adjacent cells.

The following examples are related to ManifoldCellIntGVF and ManifoldCellRealGVF. In order to make harmonic functions on manifolds, we require that the function be gradually varied. Some discussions about harmonic functions and gradually varied functions were presented in [10]. In the applications in this paper, the maximum slope of two sample points are used to determine the value levels. This treatment will guarantee that the sample point set satisfies the gradual variation condition.

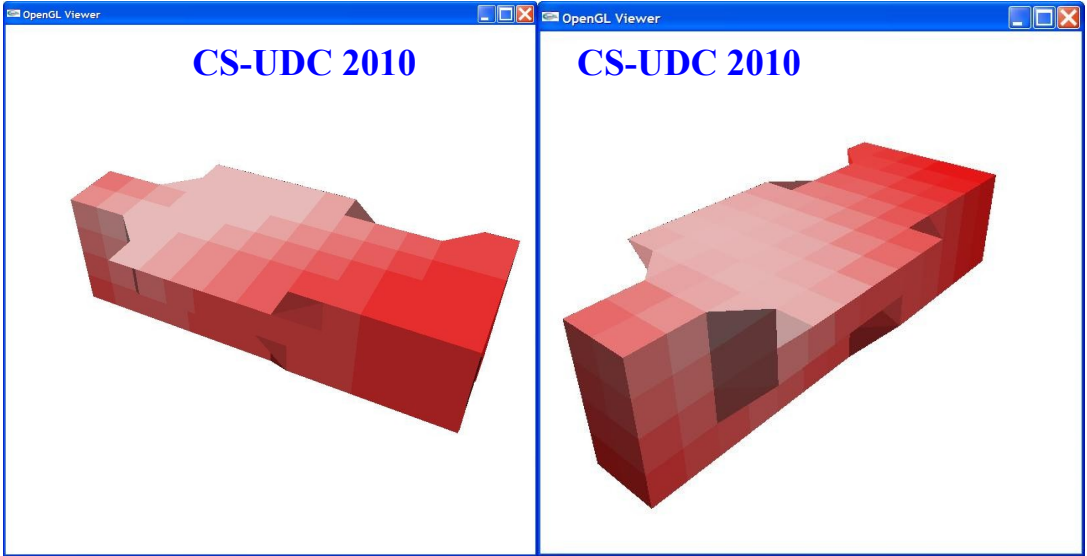

Fig. 10. Using 7 Points to fit the data on 3D surface: (a) the GVF result. (b) The Harmonic fitting (iterations few times) based on GVF.

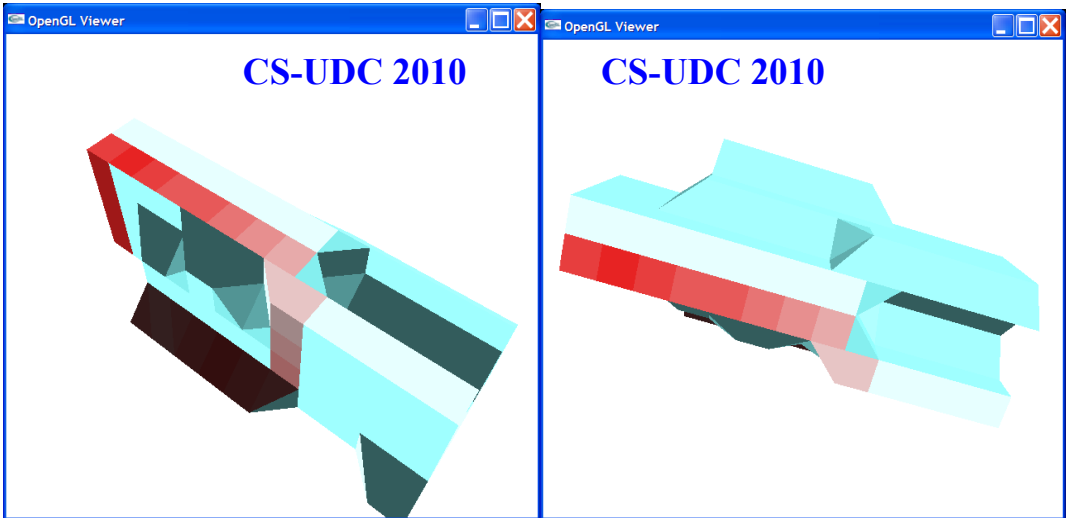

(a)                                      (b)

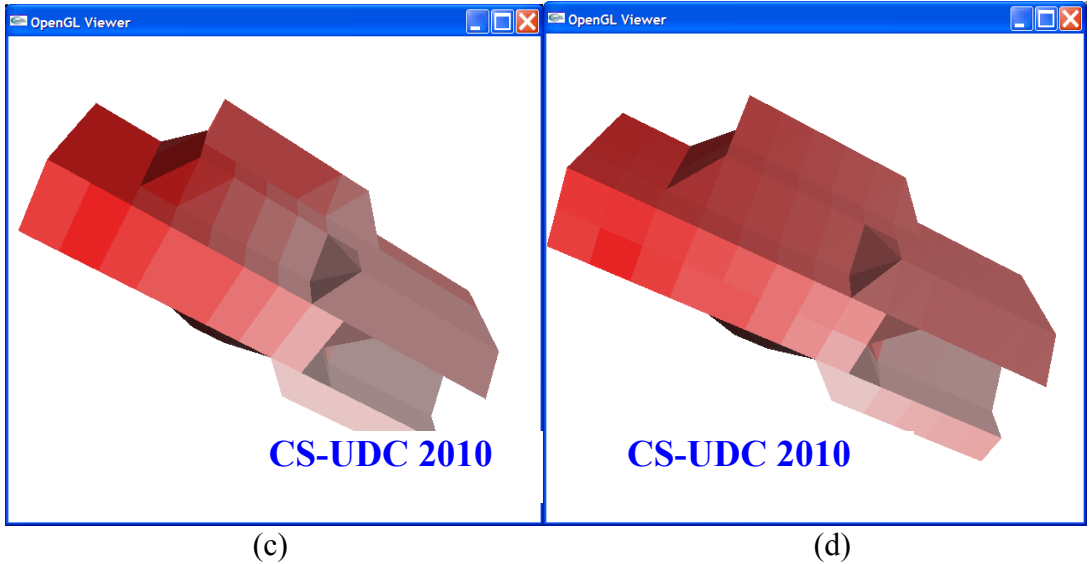

(c)  (d)

Fig. 11. The selected cells form a boundary curve that is gradually varied: (a) and (b) are two displays for the guiding points (cells). (c) The GVF result. (d) The Harmonic fitting based on GVF (100 iterations).

The smooth function on piece-wise linear manifolds is often studied in computer graphics. The main difference between our method and existing methods in graphics is apparent when we do not know the value of each vertex (or face). Our method can solve this key problem of interpolation. After the values at all vertices are known for each vertex, we can construct the smooth functions for geodesic paths that cover all edges and link to a vertex. Then we can use vector space calculus to get the unknown points. A partition of unity may be used here for getting a smooth solution on the partition edges (similar to some finite element algorithms). We can also use the existing methods for obtaining the smooth functions [11][16].

In addition, we have used our algorithm to fit the data with 7 sample points on the manifold, as described earlier in Fig. 10. We also use discrete harmonic functions to fit the data [10]. If boundary is known, then the harmonic solution will be unique. We test an example for such a case in Fig. 11.

## 5. More General Consideration

A local Lipschitz function on a manifold might not be able to be fully implemented by one set of rational or real numbers $A_1 < A_2 < \ldots < A_n$. A common method in image processing is called pyramid representation based on the intensity of an image. For a time sequence it is also well-known that a curve can be represented by a tree. It may be related to the Morse (height) function in topology. Using a tree instead of a chain for range space would solve the problem completely in theory. This method was very well prepared since in 1990, Chen proved a theorem that states that any graph can normally be immersed into a tree [6][5]. This means that there is a gradually varied extension from a graph into a tree. A more detailed description of this approach can be found in [2]. Two $A_i$ values that are the same but are located in different braches of a tree means that they represent different elements. An example is shown below:

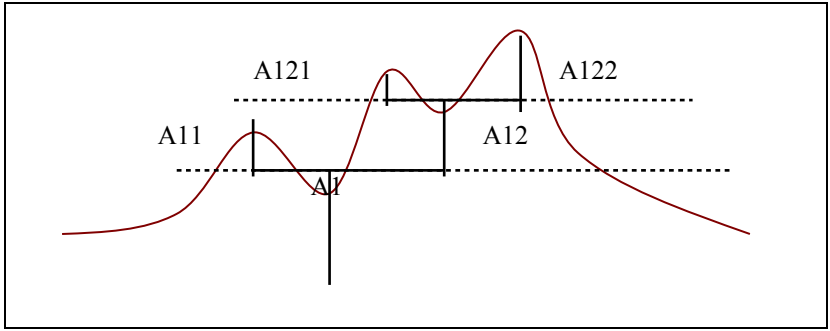

Fig. 12. A representation of range tree: Each A is a sequence of real numbers.

In Fig. 12, if A1= {B0<…<Bk}, A11={C0<…<Cr}, A12={D0<…<Dt}, then Bk=C0=D0. In this example, B1 may have the same numerical value as C3, but they represent different elements in the tree. For a real world application, since the curve is not known, a classification technique will be used for grouping the sample points (again we can use decomposition technology). It may cause some uncertainty; however, for many real world problems this method is highly applicable.

## 6. Comparison Between Gradually Varied Functions and McShane-Whitney Mid Function

The McShane-Whitney extension theorem (perhaps we should call it the McShane-Whitney-Kirszbraun Theorem since they all independently published a paper in 1934) says that a Lipschitz function f on a subset J of a connected set D in a metric space can be extended to a Lipschitz function F on D. McShane gave a constructive proof for the existence of the extension in [15]. He constructed a minimal extension (INF) that is Lipschitz. So using the same technique, it is easy for someone to construct a maximum extension (SUP). It is obvious that neither INF nor SUP can be directly used in data reconstruction. However, *F=(INF+SUP)/2* is a reasonable function that maintains good properties including each value at a point contained in the convex of guiding points. Also, *F* is a Lipschitz extension. For the purposes of this paper, let's call *F=(INF+SUP)/2* the *McShane-Whitney mid function.* The result of applying this function to our data sets is shown below. (See Fig. 13)

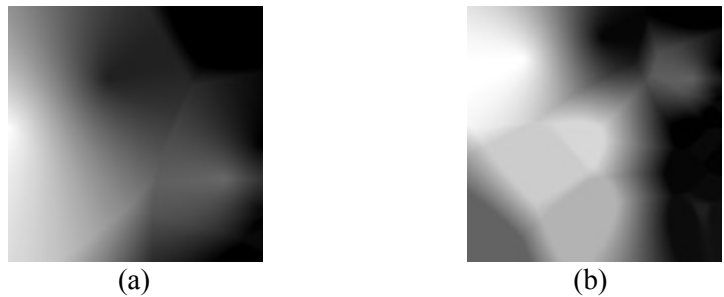

(a)          (b)

Fig. 13. McShane-Whitney mid extensions: (a) using the sample data set of Fig. 3, (b) using the sample data set of Fig. 5.

We can see that Fig. 13(a) is better than Fig. 3. However, it is not a good fit compared to Fig. 4. Fig. 13(b) is a reasonable reconstruction compared to Fig. 5. The fitting is dominated by the Lipschitz constant, so the extension function is "controlled" by sample points with bigger values.

A more detailed analysis will be given in a future paper. Fig. 13(b) looks relatively better than Fig. 13(a) since more data points (29 samples) are involved in the calcualtion. It is also more *Lipschitz*.

The difference between our algorithm and the McShane-Whitney construction method is more than the theoretical difference. Our algorithm designed based on the proof of the theorem in Section 2 is discrete and dynamic. Our method also adjusts locally. The McShane-Whitney extension method seems to be mainly for theoretical purposes. It is simple and may be fast in calculation. However, it loses flexibility especially since it is only for the Lipschitz function. In this paper, we have shown the local Lipschitz function extensions.

## 7. Summary

To get a smoothed function using gradual variation is a long time goal of our research. Some theoretical attempts have been made before, but struggled in the actual implementation. The author was invited to give a talk at the Workshop on the Whitney's Problem organized by Princeton University and the College of William and Mary in August 2009. He was somewhat inspired and encouraged by the presentations and the helpful discussions with the attendees at the workshop.

The purpose of this paper is to present some actual examples and related results using the new algorithms we designed in [2]. The author welcomes other real data sets to further examine the new algorithms. The implementation code is written in C++. Li Chen's website can be found at www.udc.edu/prof/chen.

*Acknowledgements:* This research has been partially supported by the USGS Seed Grants through the UDC Water Resources Research Institute (WRRI) and Center for Discrete Mathematics and Theoretical Computer Science (DIMACS) at Rutgers University. Professor Feng Luo suggested the direction of the relationship between harmonic functions and gradually varied functions. UDC undergraduate Travis Branham extracted the application data from the USGS database. Professor Thomas Funkhouser provided helps on the 3D data sets and OpenGL display programs. The author would also like to thank Professor C. Fefferman and Professor N. Zobin for their invitation to the Workshop on the Whitney's Problem in 2009.